\documentclass[16pt]{article}
\usepackage{amsmath,amsfonts,amssymb}
\usepackage{verbatim}
\usepackage{latexsym}
\newtheorem{theorem}{Theorem}[section]
\newtheorem{lemma}[theorem]{Lemma}
\newtheorem{proposition}[theorem]{Proposition}
\newtheorem{definition}[theorem]{Definition}

\newtheorem{remark}[theorem]{Remark}
\newtheorem{remarks}[theorem]{Remarks}

\newcommand\RR{{\Bbb R}}

\newcommand\NN{{\Bbb N}}
\newcommand\ZZ{{\Bbb Z}}

\newcommand\HH{{\Bbb H}}

\begin{document}
\title{ Mexican Hat Wavelet on the Heisenberg Group
}

\author{
 \thanks{Research supported by the German Academic  Exchange Service (DAAD).
 }
Azita Mayeli \\
\footnotesize\texttt{{mayeli@ma.tum.de}}
}

 \date{}

 \maketitle

\begin{abstract}

In this article    wavelets (admissible vectors)  on the Heisenberg group $\HH$ are studied from the point of view of \textit{Calder\'on}'s formula. We shall define  {\bf Calder\'on admissible } vectors in Definition 
\ref{eq:coroner} . Further     in Theorem  \ref{a1} we show  that   for the    class of  Schwartz functions the Calder\'on admissibility 
condition is equivalent to the  usual admissibility property  which will be  introduced in    this work. \\
  Furthermore  motivated by a well-known example on the real line, the  \textit{ Mexican-Hat} wavelet, we demonstrate the existence and construction
of an analogous  wavelet   on the Heisenberg Lie group with 2   vanishing  moments,   which  together with all of its derivatives 
has ``Gaussian" decay. The precise  proof can be found in  Theorem \ref{a6}.
\vspace{.5cm}

\footnotesize{
\begin{tabular}{lrl}
{\bf  Keywords.}  &  \multicolumn{2}{l} {\em Wavelets,  admissible vectors,  Schwartz functions,
Heisenberg groups, sub-Laplacian}\\
 &  \hspace{-.1cm}  {\em operator, Rockland operator, heat kernel.}\\
 \end{tabular}
 \vspace{.3cm}
 
 \begin{tabular}{lrl}
  &&  \hspace{-1cm}{\bf  AMS  Subject Classification (2000).}{ Primary 42C40,  42B20, 22E25}
\end{tabular}
}

     \end{abstract}

\section{Introduction and Definitions}
\label{sec:IntroductionAndPrelimianries}

   The Heisenberg group $\HH$ is  a Lie group with underlying manifold $\RR^3$. We denote points in $\HH$ by $(p,q,t)$ with $p,q,t\in \RR$, and define the group operation by 
\begin{align}\label{group-action}
(p_1,q_1,t_1)\ast (p_2,q_2,t_2)= \big(p_1+p_2,\; q_1+q_2,\; t_1+t_2+\frac{1}{2}(p_1q_2-q_1p_2)\big).
\end{align}
It is straightforward to verify that $\ast$ is a group operation.
We can identify both $\HH$ and its Lie algebra $\mathfrak{h}$ with $\RR^3$, with group operation given by (\ref{group-action}) and Lie bracket given by 
\begin{align}\notag
\left[(p_1,q_1,t_1),(p_2,q_2,t_2)  \right]= (0,0,p_1q_2-q_1p_2).
\end{align}
The Haar measure on the Heisenberg group $\HH= \RR^3$ is the usual Lebesgue measure.
More precisely, the Lie algebra $\mathfrak{h}$ of the Heisenberg group $\HH$ has a basis $\{X,Y,T\}$, which we may think of as left invariant differential operators on $\HH$;  where $[X,Y]=T$  and  all other  brackets are zero,  and where  the exponential function $\exp:\;\mathfrak{h}\rightarrow \HH$ is the identity, i.e., 
\begin{align}\notag
\exp\;(pX+qY+tT)= (p,q,t).
\end{align}
The  action of  $\mathfrak{h}$ on space $C^\infty(\HH)$ via   the 
  left invariant differential operators  $\{X,Y,T\}$  is defined by the following formula:\\
   Suppose $f\in C^\infty(\HH)$, then  
   \begin{align}\label{XYZ}
(Xf)(p,q,t)=&\frac{d}{dp}f(p,q,t)-\frac{1}{2}q\frac{d}{dt}f(p,q,t),\\\notag
   (Yf)(p,q,t)=&\frac{d}{dq}f(p,q,t)+\frac{1}{2}p\frac{d}{dt}f(p,q,t),\\\notag
   (Tf)(p,q,t)=& \frac{d}{dt}f(p,q,t).
   \end{align}
  Our definition of continuous wavelet transform for the Heisenberg group 
  will be from  the representation  point of view adapted  from the case $\RR$.     
       For the construction of wavelet transform  one needs  a one-parameter group  of dilations for $\HH$. Here we consider  $H:=(0,\infty)$ as  the  
         one-parameter dilation group of $\HH$  which is defined as follows:
   Suppose $a>0$. Then the operator $\delta_a$ 
   defines an automorphism of $\HH$ by 
   \begin{align}\label{automorphism}
   \delta_a(p,q,t)= (ap,aq,a^2t)\hspace{.5in}\forall (p,q,t)\in\HH.
   \end{align}
    The set  $\{\delta_a:\;a>0\}$ forms a     group of automorphisms of $\HH$, called the   dilation  group for $\HH$ (for more details about such dilation groups  see for example \cite{FollandStein82}).
We denote the operation of $\delta_a$  by $\delta_a(\omega)= a\omega$ for any $\omega\in \HH$. From now on, 
   $a>0$ refers to the automorphism $\delta_a$ 
and   $H=(0,\infty)$  denotes the closed  subgroup of automorphisms of $\HH$ with operation as  in (\ref{automorphism}).\\   
  The group $H=(0,\infty)$  operates continuously  by topological automorphisms on the locally compact group $\HH$. So we can define the semidirect product $G:=\HH\rtimes (0,\infty)$, which is a locally compact topological  group with  the product topology. Elements of $G$ can be written as $(\omega, a)\in \HH\times (0,\infty)$ and the group operation on $G$  is defined by 
  \begin{align}\notag
  (\omega,a)(\acute{\omega},\acute{a})= (\omega(a\acute{\omega}), a\acute{a})\quad \forall \omega, \acute{\omega}\in \HH\; \;\text{and} \;\;\forall a, \acute{a}>0.
  \end{align}
  $G$ is a non-unimodular group and its left Haar measure is given by $d\mu_G(\omega,a)= a^{-5}d\omega da$. 
  Analogously  to the  situation on $\RR$, for $a>0$ the dilation operator $D_a$ is defined by
   $D_af(.)=  a^{-2} f(a^{-1}.)$ and for $\omega\in \HH$, $L_\omega $ denotes the   left translation where $L_\omega f(.)= f(\omega^{-1}.)$ for any $f$ defined on $\HH$ and $\omega\in \HH$. 
 \begin{definition}\label{quasi-regular representation}
 For any $(\omega,a)\in G$ and $f\in L^2(\HH)$ define 
 \begin{align}\label{quasi-rep}
 (\pi(\omega,a)f)(\upsilon):=L_\omega D_af(\upsilon)= a^{-2}  f(a^{-1}(\omega^{-1}.\upsilon)).
 \end{align}
 It is easy to prove that $\pi$ is a strongly  continuous unitary representation of $G$ which  acts  on $L^2(\HH)$ by  $(\ref{quasi-rep})$. This representation is 
 called the  ``quasi-regular representation".
 \end{definition}
  Next we give the definition of admissible vectors in $L^2(\HH)$, which arises from the action  of the quasi-regular representation on  $L^2(\HH)$ by  (\ref{quasi-rep}).
 \begin{definition}\label{admissible-vector}
 For any $\phi\in L^2(\HH)$ the associated  coefficient operator $V_\phi$ is defined  on $L^2(\HH)$ by 
 \begin{align}\notag
  V_\phi (f)(\omega,a)= \langle f, \pi(\omega,a)\phi\rangle \quad \forall f\in  L^2(\HH), \; (\omega, a)\in G.
 \end{align}
 $\phi$ is called {\bf admissible} if $V_\phi$ maps  $L^2(\HH)$ into  $L^2(G) $   isometrically    up to a constant, i.e., 
  \begin{align}\label{isometric}
  \parallel f\parallel^2= const. \int_\HH\int_0^\infty \mid V_\phi(f)(\omega,a)\mid^2 a^{-5}dad\omega \quad \forall\; f\in L^2(\HH),
  \end{align}
  where the constant is positive   and only depends on $\phi$. Then $V_\phi$ is called  a
  {continuous wavelet transform} and $V_\phi(f)$ is called  the  continuous wavelet transform of   $f$. 
 \end{definition}
 One of the important consequence of the isometry given by  formula (\ref{isometric}) is that a function can be reconstructed from its wavelet transform by means of the \textit{``resolution identity"}, i.e, formula (\ref{isometric}) can be read as 
  \begin{align}
  f= const. \int_\HH\int_0^\infty \langle f, \pi(\omega,a)\phi\rangle \pi(\omega,a)\phi\; a^{-5}dad\omega\quad \forall\; f\in L^2(\HH),
  \end{align}
  which the  convergence of the integral is understood in the weak sense.\\
  
 The most  importance of wavelet  theory is  its
  microscope effect, i.e,  by choosing a suitable wavelet $\phi$, as the lens, one can obtain     information about the local regularity of  argument functions $f$ in $L^2(\HH)$. This information   is  obtained  from the wavelet  coefficients    $\langle f, \pi(\omega,a)\phi\rangle$ when for instance these  coefficients have a fast decay when $a\rightarrow 0$. 
  Note that in the following definition we take the Schwartz functions on the Heisenberg group to be  the Schwartz functions on  $\RR^3$. 

 The existence of admissible vectors for the quasi-regular representation of $G:=N\rtimes H$ on $L^2(N)$ is already proved in F\"uhr's book \cite{Fuehr05},  where $N$ is a homogeneous Lie group and $H$ is  a one-parameter group  of dilations for $N$ (for the definition of homogeneous groups see for example \cite{FollandStein82}). However, the existence of  smooth fast-decaying wavelets was left open.\\
 Our work establishes existence of admissible radial Schwartz vectors for  the case $N=\HH$ and $H=(0,\infty)$. 

The existence of admissible vectors in  closed subspaces of  $L^2(\HH)$ was studied   in \cite{LiuPeng97}. The authors consider the unitary reducible representation $U$ of a non-unimodular group $P$ on $L^2(\HH)$. They  
  decompose $L^2(\HH)$ into an infinite  direct   sum of the irreducible  invariant closed subspaces, $\mathcal{M}_n$, under  the  representation  $U$ on $L^2(\HH)$. Then they  show that the restriction of $U$ to  these subspaces is square-integrable, i.e, each subspace $\mathcal{M}_n$ contains at last one nonzero  wavelet vector with respect to $U$. Furthermore  the authors give a characterization of the admissibility  condition in the  irreducible  invariant closed subspaces $\mathcal{M}_n$ in the terms of the  Fourier transform. 
  But it seems that  it was not trivial for the authors  to show the existence of an admissible vector for all of  $L^2(\HH)$. In contrast, as we will seen soon,  our  work  first provides a characterization  of 
     admissible functions  in the Schwartz space on the Heisenberg group,  and then presents an explicit example  of  an  admissible function.\\

It  is particularly remarkable  that the   representation $U$ in \cite{LiuPeng97} is unitarily 
equivalent to the  direct sum of irreducible representations, which   are  all square integrable. Hence by 
      Corollary 4.27 in \cite{Fuehr05}, the representation $U$ on $L^2(\HH)$ is square integrable, i.e, there exists an admissible vector in $L^2(\HH)$. Therefore by  relying on this consequence of \cite{Fuehr05}
      we are aware of existence of admissible functions in $L^2(\HH)$. In this work we want to study the admissibility  condition for  functions in $L^2(\HH)$ with respect to the quasi-regular representation, and obtain a  concrete example of a Schwartz wavelet with some nice additional  properties.\\
      
     As we saw above,  it seems that   the  study  of     wavelets on the Heisenberg group from
     the  representation theory
      point of view  is  the  usual approach to the subject. 
      Since that method does not easily provide an explicit example for  the group, we    look at  wavelets  through an equivalent  approach, which we will discuss  in the next section. \\
      
       We organized the  new results contained in the work  as follows:   In section 
       \ref{CalderonAdmissibleFunctions}  we define the  ``Calder\'on admissibility'' of  a function on 
       the Heisenberg group   and then in Theorem \ref{a1} we show that the (accepted)  definition of  admissibility  is consistent with our usage of the word of wavelets as in   Definition (\ref{admissible-vector}).  This theorem 
   provides  the characterization of     wavelets in the class of  Schwartz functions on the Heisenberg group.            In section \ref{sec:TheExample} we construct an explicit example of  a Schwartz wavelet  with   two    vanishing moments, such that it and all of its derivatives have ``Gaussian'' decay.  (We say  a function $F$ on $\HH$ has  ``Gaussian''  decay if for some $C, c>0$, 
   \begin{align}\notag
   \mid F(p,q, t)\mid \leq C \exp c(p^4+q^4+t^2)^{-1/2}\quad \forall (p,q,t)\in \HH.)
   \end{align}

  \section{Calder\'on Admissible Functions }     
  \label{CalderonAdmissibleFunctions} 
 As mentioned before, the existence of an admissible vector for $L^2(N)$ is proved in \cite{Fuehr05}, where $N$ is a homogeneous group, for the  quasi-regular representation of $G:=N\rtimes H$ on $L^2(N)$. Here  $H$ is a one-parameter group of dilations  of $N$.
The existence of such vectors for the case $N:=\RR^k$ and $H < GL(k,\RR)$ has  recently been studied by different authors,  for instance for $k\in \NN$ and $H$ as a closed subgroup of $GL(k,\RR)$ by Hartmut F\"uhr  in \cite{Fuehr96} and \cite{Fuehr98},   and  for the case $k=1$ and $H:=\ZZ$ by  the authors in  \cite{MallatZhong92}. The case $N:=\HH$  and  $H:=\RR$ as a one-parameter group of dilation is considered by  \cite{LiuPeng97}. 
In this section we prove the existence of  Schwartz admissible vectors for the case $N:=\HH$ and $H=(0,\infty)$ by applying the following definition:
\begin{definition}\label{eq:coroner} Let $ \phi\in \mathcal{S}(\HH)$ and $\int \phi =0$. Then $\phi$ is called 
 \textit{\bf Calder\'on} admissible if for any \;$ 0<\varepsilon < A$ and $g\in \mathcal{S}(\HH)$   
\begin{align}\label{admiss-vec-for-stratified}
g\ast \int_\varepsilon^{A}  \tilde{\phi}_a\ast \phi_a   \; a^{-1}da  \rightarrow cg\quad \text{as} \;\;\varepsilon\rightarrow 0;\; A\rightarrow \infty \quad 
\end{align} 
holds 
in the sense of tempered distributions
 where 
  $c$ is a nonzero constant and  $\phi_a(\omega)= a^{-4}\phi(a^{-1}\omega)$.   
\end{definition} 
 
  In Lemma \ref{a1} below, we show  that on the Schwartz space  the   definition of admissibility  in (\ref{eq:coroner})    is equivalent to    the word 
  \textit{admissible}  in the sense of  Definition  (\ref{admissible-vector}):
    \begin{theorem}\label{a1}
 Let $\phi\in \mathcal{S}(\HH)$ and $\int \phi=0$, then $\phi$ is admissible   if and only if     $\phi$  is  \textit{Calder\'on} admissible.
 \end{theorem}
 {\bf proof:\;} Suppose  $\phi\in \mathcal{S}(\HH)$ and $g\in \mathcal{S}(\HH) $. Then according to Definition \ref{admissible-vector} we have:
\begin{align}\notag
\|V_{ {\phi}}g\|_2^2&=\int_0^\infty \int_{\HH}|\langle g,
\lambda(b) D_{a } {\phi}\rangle|^2db a^{-5 }da\\ \notag
&=\int_0^\infty \int_{\HH}|g\ast D_{a }\widetilde{\phi}(b)|^2db a^{-5 }da\\\notag
 &=\int_0^\infty \|g\ast D_{a }\widetilde{\phi} \|_{L^2(\HH)}^2a^{-5 }da\\\notag
&=\underset{\varepsilon \rightarrow 0,\;A \rightarrow
\infty}{\lim}  \int_{\varepsilon }^A  
\|g\ast D_{a }\widetilde{\phi} \|_{L^2(\HH)}^2a^{-5 }da\\\notag
  &=\underset{\varepsilon \rightarrow 0,\;A \rightarrow
\infty}{\lim}  \int_{\varepsilon }^A \langle g\ast D_{a } \widetilde\phi\;, \;g\ast
D_{a }\widetilde\phi\rangle a^{-5 }da\\ \notag
 &=\underset{\varepsilon \rightarrow 0,\;A \rightarrow
\infty}{\lim}  \int_{\varepsilon }^A  \langle g \;,\;g \ast D_{a } \widetilde{\phi} \ast {
D_{a }{\phi}}\rangle a^{-5 }da\\ \notag
    &= \underset{\varepsilon \rightarrow 0,\;A \rightarrow
 \infty}{\lim}\langle g\;,\;
  g\ast \int_{\varepsilon }^A   D_{a } \widetilde{\phi} \ast
 D_{a } {\phi}a^{-5}da\rangle \\\label{the-last-one}
 &= \underset{\varepsilon \rightarrow 0,\;A \rightarrow
 \infty}{\lim}\langle g\;,\; g\ast \int_{\varepsilon }^A    \widetilde{\phi_a} \ast
  {\phi}_a a^{-1}da\rangle.   
\end{align}
If $K_{\varepsilon, A}= \int _\varepsilon^A  \widetilde{\phi_a} \ast
  {\phi}_a a^{-1}da$ then $K= \underset{\varepsilon \rightarrow 0,\;A \rightarrow
 \infty}{\lim}K_{\varepsilon, A}$ exists in $\mathcal{S}' (\HH)$, $C^\infty$ away from $0$  and is homogeneous of degree $-4$, by Theorem 1.65 in \cite{FollandStein82}.
Therefore  if  $g\in \mathcal{S}$, $g\ast K_{\varepsilon, A}\rightarrow g\ast K$ pointwise and for some $N, C$
\begin{align}\notag
\mid (g \ast K_{\epsilon,A})(x)\mid 
\leq C(1+\mid x\mid)^N \quad \quad \text{for all } x,\varepsilon,
A.
\end{align}
 
 Using the dominated convergence theorem in (\ref{the-last-one}), if $g\in \mathcal{S}(\HH)$, then   
 \begin{align}\label{inequality}
 \langle g\;,\;\underset{\varepsilon \rightarrow 0,\;A \rightarrow
\infty}{\lim}g\ast \int_{\varepsilon }^A    \widetilde{\phi_a} \ast
 {\phi}_a a^{-1}da\rangle = \langle g, g\ast K\rangle\leq C\parallel g\parallel_{2}^2
 \end{align}
 since the map $g\rightarrow g\ast K$ is bounded on $L^2(\HH)$. Thus $V_\phi$ maps
$\mathcal{S}(\HH)$
to $L^2(G)$ and has a unique bounded extension to a map from $L^2(\HH)$ to $L^2(G)$. But if
$g_k\rightarrow g$ in
$L^2(\HH)$, surely $V_\phi g_k\rightarrow V_\phi g$ pointwise , so this extension can be none
other than $V_\phi$.  Accordingly (\ref{inequality}) holds for all $g\in L^2(\HH)$. We thus have

\begin{align}\notag
\parallel V_\phi g\parallel_{2}= \parallel g\parallel_{2}\quad  \forall g\in L^2
&\Longleftrightarrow \langle g ,g\ast K\rangle  = \langle g,g\rangle  \quad \forall g\in L^2 \\\notag
&\Longleftrightarrow g\ast K=g \quad \forall g\in L^2 \\\notag
&\Longleftrightarrow K=\delta \quad \text{up to a constant.}\quad
\end{align}

as desired.  (In the second implication, we have used polarization.)  
This completes the proof. \quad$\Box$ \\

In  the next Proposition we will  obtain a sufficient condition for Schwartz functions to be admissible which is one of the chief tools for the proof of our  main theorem.  
 \begin{proposition}\label{a3}
Suppose $\phi,\psi \in \mathcal{S}(\HH)$, so that $\int \phi=0$ and   $\int\psi\not=0,$
and for some constants $k, c > 0$ and non-zero real number $q$   one has
  $\widetilde{\phi}_{a^q} \ast  {\phi_{a^q}}=-ac\frac{d}{da} \psi_{ka^q}$. Then $\phi$ is
 admissible.
 \end{proposition}
 {\bf proof:\;} Suppose $g\in \mathcal{S}(\HH)$ and $0<\varepsilon <A <\infty$. By  changing the  coordinate $a$ to $a^q$  and  by using the assumption that $\widetilde{\phi}_{a^q} \ast  {\phi_{a^q}}=-ac\frac{d}{da} \psi_{ka^q}$, 
  we  can write 
 \begin{align}\notag
g\ast \int_{\varepsilon }^A    \widetilde \phi_a \ast
\phi_a\;a^{-1}da&= q \;g\ast \int_{\varepsilon^{1/q}}^{A^{1/q}} \widetilde \phi_{a^q}\ast  \phi_{a^q}
\;a^{-1}da\\  \label{a4}
&= (qc) \;g\ast \int_{\varepsilon^{1/q}}^{A^{1/q}}
\left(-a\frac{d}{da}\right)\psi_{ka^q}\;a^{-1}da\\ \notag  
&=( qc) \;g\ast \int_{\varepsilon^{1/q}}^{A^{1/q}}\left(- \frac{d}{da}\psi_{ka^q}\right)
da\\\notag  
&=(-qc)\left(g\ast (  \psi_{ A k}  - \psi_{\varepsilon k})\right)\\\notag  
&=(qc) \left(g\ast (\psi_{\varepsilon k}-\psi_{ A k} )\right)
\end{align}
Since $\int\psi\not=0$, then from Proposition 1.20 \cite{FollandStein82} we have:
\begin{align}\label{a41}
\underset{\varepsilon \rightarrow 0}{\lim}\;g\ast\psi_{\varepsilon k}=g\int \psi,\hspace{.2in} \text{in} \; L^2\text{-norm}.
\end{align}
On the other hand one  can write:
\begin{align}\label{a5}
\|g\ast \psi_{Ak}\|_2&\leq \|g\|_1\|\psi_{kA}\|_2=(kA)^{-\frac{1}{2}}\|g\|_1\|\psi\|_2
\end{align}
which shows $  g\ast \psi_{Ak}\rightarrow 0$    in $L^2$-norm as $A\rightarrow \infty$.
 Now  applying  (\ref{a41}) and   (\ref{a5}) in  (\ref{a4})  
 \begin{align}\notag
g\ast \int_{\varepsilon }^A    \phi_a \ast 
\tilde{\phi}_aa^{-1}da \rightarrow (qc) g\int \psi \quad as \quad \varepsilon\rightarrow 0,\; A\rightarrow \infty \quad  \text{in} \;L^2\text{-norm},
\end{align}
as   desired.\quad $\Box$
\vspace{1cm}
  
  To present  our next    main result,  we have to recall
  some basic   definitions  first:  \\ 
Suppose $L= -(X ^2+Y^2)$ is the 
 sub-Laplacian operator, where $X $ and $Y$ are the 
left-invariant  vector fields on the Heisenberg  group which have been defined earlier  in (\ref{XYZ}). The \index{operator!heat kernel}heat kernel operator associated to $L$ is the differential operator $ \frac{d}{dt} +L$ on $\mathcal{\HH}\times \RR$, where $\frac{d}{dt}$ is the coordinate vector field on $\RR$ (one can consider this coordinate as the time coordinate). For  the heat operator we recall here    Proposition 1.68 of \cite{FollandStein82} for the Heisenberg group.
\begin{proposition}\label{heat-kernel}
There exists a unique   $C^\infty$  function $h$ on $\HH\times (0,\infty)$,   for which   the following properties hold :
\begin{enumerate} 
	\item $(\frac{d}{dt}+L)h=0 \; \;\text{on } \;\HH\times (0,\infty)$
	\item $ h(\omega,t)\geq 0, \; h(\omega,t)= h(\omega^{-1},t)\quad \forall(\omega,t)\in \HH\times (0, \infty)$
		\item $ \int h(\omega,t)d\omega=1 \; \text{for} \;  t>0  $
	\item $ h(.,s)\ast h(.,t)=h(.,s+t) \quad \forall \; s,t>0$
	\item  $r^{4}h(r\omega,r^2t)=h(\omega,t)\quad \forall \;\omega\in \HH,\; t,r>0,$ \\ 
	$($note that  here $r\omega$   is understood as the result of applying the 
	   automorphism $\delta_r$ to $\omega$,  that is,  $\delta_r(\omega))$. 
\end{enumerate} 
The solution $h$ is known as the   
\textit{\bf  heat kernel}.
\end{proposition} 
\begin{remarks}

\begin{itemize}
\item Note that the   proposition  \ref{heat-kernel} has been  proved for the stratified groups in \cite{FollandStein82}. 
\item
 Here the interval $(0,\infty)$ has nothing to do with the  one-parameter group of dilations which has   been  introduced earlier. 
  One may consider  it as  a time interval.
  \end{itemize}
  \end{remarks}
 The idea of this section is  to 
  apply Proposition $\ref{a3} $  to  $\phi(x)=L h(x,1)$  to show that the function $\phi$ is an admissible vector. For that reason   here we first need   to compute the  dilates  of functions $h(.,1)$ and  $Lh(.,1)$. 
 \begin{lemma}
 For any $a>0$ and $\omega\in \HH$ we have  
 \begin{center}
 $ h(\omega,1)_a=a^2h(\omega, a^2)$\; and \;
  $Lh(\omega,1)_a= a^2Lh(\omega,a^2)$.
  \end{center}
 \end{lemma}
 {\bf Proof:\;} Suppose $a>0$ and $\omega\in \HH$. Applying    $\#5$ in Proposition  \ref{heat-kernel}  one finds: 
   \begin{align}\label{dilation-of-h}
  h(\omega,1)_a= a^{-4} h(a^{-1}\omega,1)=a^2h(\omega, a^2).
\end{align}
 Similarly  by applying  $\#1$ and $\#5$   in Proposition  \ref{heat-kernel} for $Lh(.,1)$ we have:
\begin{align}\label{dilation-of-Lh}
Lh(\omega,1)_a&= a^{-4}Lh(a^{-1}\omega,1)\\ \notag 
&= -a^{-4}\frac{d}{dt}h(a^{-1}\omega,t)\big|_{t=1}\\ \notag
&=-\frac{d}{dt}h( \omega,a^2t)\big|_{t=1}\\ \notag
&= a^2Lh(\omega,a^2).\quad \Box
\end{align}
 
   \section{Mexican Hat  Wavelet on $\HH$}
\label{sec:TheExample}
 The purpose  of this section is  to show the   \textit{Calder\'on} admissibility   of  the  Schwartz  function $\phi= Lh(.,1)$, which   will be  stated   in    Theorem \ref{a6} as the other  main result of this work.   But first we make  the following remark:
   \begin{remark}
  On the real line $\RR$, the \textit{ heat kernel}, as our motivating  example,   is given by $h(x,t)= \frac{1}{\sqrt{4\pi}} e^{-\frac{x^2}{4t}}$,  and   $h(x,1)= \frac{1}{\sqrt{4\pi}} e^{-\frac{x^2}{4}}.$ The second derivative of the 
  Gaussian  is an often employed wavelet, 
     the  \textit{Mexican-Hat} wavelet. This function (and its  dilated and translated copies) has a shape similar to a Mexican hat (for more see for example \cite{Daube92})
     This wavelet has two vanishing moments  and evidently   it and all its derivatives have Gaussian decay.  Our goal in the next   theorem is  prove the existence of a Mexican-Hat wavelet on the Heisenberg group with   similar   properties.
  \end{remark}
   
\begin{theorem}\label{a6} The 
Schwartz  function $\phi(\omega)= Lh(\omega, 1)$ is admissible and  it and all of its derivatives have ``Gaussian" decay.
\end{theorem}
 {\bf Proof:\;}  First we shall  show that $\tilde{\phi}=\phi$. This  is easy to see 
since  for any  $\omega\in\HH$ and $t>0$ 
by applying   $\#2$ in Proposition 
 \ref{heat-kernel},  we find
 \begin{align}\notag
\widetilde{Lh}(\omega,t)=
  -\frac{d}{dt} \tilde{h}(\omega,t) 
=-\frac{d}{dt}h(\omega,t)=Lh(\omega,t). 
\end{align}

To prove  the theorem,  it is sufficient to show that for  the function $\psi= h(\omega,1)+Lh(\omega,1)$ the relation
\begin{align}\label{a8.1}
  \widetilde\phi_{\sqrt{a}}\ast \phi_{\sqrt{a}}=\phi_{\sqrt{a}}\ast \phi_{\sqrt{a}}=-ca \frac{d}{da}\psi_{\sqrt{2a}} 
  \end{align}
  holds.
 Hence
  by applying Proposition  \ref{a3} we will   get our assertion. \\
     Using the relations  (\ref{dilation-of-h}) and  (\ref{dilation-of-Lh}),  we find 
\begin{align}\label{a9}
\phi_{\sqrt{a}}\ast\phi_{\sqrt{a}}&=(Lh(.,1))_{\sqrt{a} }\ast(Lh(.,1))_{\sqrt{a}}\\
\notag &=aLh(.,a)  \ast aLh(.,a).
\end{align}
But for any $a,b > 0$, 
\begin{align}\notag
Lh(.,a)*Lh(.,b)& = \frac{d}{da}h(.,a)* \frac{d}{db}h(.,b)\\\notag
&= \frac{d}{da} \frac{d}{db} h(.,a+b) \\\notag
&= L^2 h(.,a+b).
\end{align}
 If $a=b$
we get   
\begin{align}
\phi_{\sqrt{a}}\ast\phi_{\sqrt{a}}=a^2 L^2 h(.,2a),
\end{align}

while, by using (\ref{dilation-of-Lh}), we find  
\begin{align}\notag
\psi_{\sqrt{2a}}&=(h(.,1))_{\sqrt{2a}}+(Lh(.,1))_{\sqrt{2a}}\\ \notag 
&=h(.,2a)+2a Lh(.,2a). 
\end{align}
Observe that the  derivative of $\psi_{\sqrt{2a}}$ with respect to the parameter $a$   is computed  as follows :
\begin{align}\label{a11}
\frac{d}{da}\psi_{\sqrt{2a}}&=\frac{d}{da}h(.,2a)+2Lh(.,2a)+2a\frac{d}{da}Lh(.,2a)\\ \notag 
&=2\frac{d}{d2a}h(.,2a)+2Lh(.,2a)+4a\frac{d}{d2a}Lh(.,2a)\\ \notag
 &=-4aL^2h(.,2a). 
\end{align}
Comparing the  equations  (\ref{a9}) and  (\ref{a11}), we see that  the relation  (\ref{a8.1})  holds for
$\phi$,   $\psi$,  and for  $c=4$, as desired.\\
The fact that the function $\phi$ has the 
    property  that  it and all of its derivatives have   ``Gaussian'' decay is known by the work of 
Jersion and Sanchez-Calle \cite{JerSan86} and of Varopoulos \cite{Varopo88} and hence we are done.\quad $\Box$
       
   \section{ Some Remarks}
   \label{remarks}
   \begin{enumerate}
     \item In this article, we  provided  our  results for  the  Heisenberg group, $\HH\simeq \RR^3$,    only for the sake of simplicity;   evidently 
     our main results in this work hold for the Heisenberg group  $\HH^n$  for  any   $n$ also, i.e.  with the underline manifold  $\RR^{2n+1}$.
     \item   Observe that the main results of this article   can   also be achieved for  the general case of    stratified  Lie groups of any homogeneous degree, since Proposition
     \ref{heat-kernel} has been stated 
  for  this  class of  groups. (For more about  stratified groups and   homogeneous degree, we refer the interested  reader to  \cite{FollandStein82}.)
      \item  Again for simplicity,  here we     considered   only   the sub-Laplacian operator  $L$   on the group, but  certainly one can see that  the results hold for any positive Rockland operator,  which is defined   for example  in 
      \cite{FollandStein82}.
  \end{enumerate}

The author would like to thank  Daryl N.Geller, G\"unter Schlichting, and Hartmut F\"uhr  for   helpful discussions.


\begin{thebibliography}{99} 
\bibitem{Christensen03} O.  Christensen, \textit{An Introduction to  Frames and Riesz Bases}, Brikha\"user, 2003. 
\bibitem{Daube92} I. Daubechies, \textit{Ten Lectures on Wavelets}, Philadelphia, Pennsylvania, 1992. 
\bibitem{FollandStein82} G.B.Folland, E.M.Stein,\textit{ Hardy spaces 
     on homogeneous groups}, Princeton University Press, Princeton, 1982.
\bibitem{Folland95} G.B.Folland,     \textit{A Course in Abstract Harmonic Analysis}, CRC Press, Boca Raton, 1995. 
     \bibitem{Fuehr96} H. F\"uhr, \textit{Wavelet Frames and Admissibility in Higher Dimensions}, J. Math. Phys. {\bf 37} (1996), 6353--6366.
        \bibitem{Fuehr98} H. F\"uhr, \textit{Continuous Wavelet Transforms with Abelian Dilation Groups}, J. Math. Phys. {\bf 39} (1998), 3974--3986. 
 \bibitem{Fuehr05}H. F\"uhr, \textit{Abstract Harmonic Analysis of Continuous Wavelet Transforms}, 
Lecture Notes in 
 Mathematics {\bf  1863},  Springer Verlag, Berlin, 2005.
 \bibitem{gm1} D. Geller and A.  Mayeli, \textit{Continuous wavelets and frames on stratified Lie
groups I}, to appear in Journal of Fourier Analysis and Applications, 2006.
 \bibitem{JerSan86} D. Jerison and A.  Sanchez-Calle, \textit{Estimates for the heat kernel
for a sum of squares of vector fields}, Indiana Univ. Math. J. 35 (1986), 835-854.
  \bibitem{LiuPeng97}{ H. Liu and L. Peng,  {\em Admissible wavelets associated
 with the Heisenberg group,} Pac. J. Math. {\bf 180} (1997), 101-123.}
   \bibitem{MallatZhong92}{S. Mallat and S.  Zhong, {\em Wavelet transform maxima and
 multiscale edges,} in Wavelets and Their Applications, M.B. Ruskai,
 G. Beylkin, R. Coifman, I. Daubechies, S. Mallat, Y. Meyer, and L. Raphael,
 eds., Jones and Bartlett, Boston, 1992, 67-104.}
 \bibitem{Mayelithesis05} A. Mayeli, \textit{Discrete and continuous wavelet 
transformation on the Heisenberg group}, Ph.D thesis, Technische Universit\"at M\"unchen, 2005.
\bibitem{Varopo88} N.  Varopoulos, \textit{Analysis on Lie Groups}, 
J. Func. Anal. 76 (1988), 346-410.
\end{thebibliography}
\end{document}